\documentclass[letterpaper, 10 pt, conference]{ieeeconf}  
\IEEEoverridecommandlockouts   
\usepackage{amsmath, amssymb}
\usepackage{algorithm, algpseudocode}												%
\usepackage[T1]{fontenc}
\usepackage{subfig}
\usepackage{stmaryrd}
\usepackage[latin1]{inputenc}
\usepackage{graphicx}
\usepackage{amssymb}
\usepackage{amsmath}
\usepackage{dcolumn}
\usepackage{bm}
\usepackage{color}
\usepackage{nomencl}

\newcommand{\vt}{\boldsymbol} 
\newcommand{\G}{\mathcal G}
\newcommand{\F}{\mathbb F}
\newcommand{\R}{\mathbb R}
\newcommand{\ones}{\boldsymbol 1}
\DeclareMathOperator{\prob}{prob}
\DeclareMathOperator{\rank}{rank}
\DeclareMathOperator{\diag}{diag}
\DeclareMathOperator{\round}{round}
\newtheorem{proposition}{Proposition}
\newtheorem{theorem}{Theorem}
\newtheorem{lemma}{Lemma}

\title{\LARGE \bf Distributed privacy-preserving network size computation:\\ 
A system-identification based method}

\author{Federica Garin and Ye Yuan%
\thanks{This is a pre-print of a paper  which will appear in the proceedings of the 52nd  IEEE Conference on Decision and Control (CDC 2013), Dec.~2013, Firenze, Italy. The appendix presented here complements the conference paper with proofs which were omitted in the proceedings due to page-number limitation.}
\thanks{This collaborative work is an outgrowth of the LCCC focus period on
Information and Control in Networks organized at Lund University in
October 2012.}%
\thanks{F.~Garin is with NeCS team, INRIA Grenoble -- Rh\^one-Alpes, France
        {\tt\small federica.garin@inria.fr}}%
\thanks{Y.~Yuan is with Department of Engineering, University of Cambridge, UK
        {\tt\small yy311@cam.ac.uk}}%
}

\begin{document}

\maketitle
\thispagestyle{empty}
\pagestyle{empty}

\begin{abstract}
In this study, we propose an algorithm for computing the network size of communicating agents. The algorithm is distributed: a) it does not require a leader selection; b) it only requires local exchange of information, and; c) its design can be implemented using local information only, without any global information about the network. It is privacy-preserving, namely it does not require to propagate identifying labels. This algorithm is based on system identification, and more precisely on the identification of the order of a suitably-constructed discrete-time linear time-invariant system over some finite field. We provide a probabilistic guarantee for any randomly picked node to correctly compute the number of nodes in the network. Moreover, numerical implementation has been taken into account to make the algorithm applicable to  networks of hundreds of nodes, and therefore make the algorithm applicable in real-world sensor or robotic networks. We finally illustrate our results in simulation and conclude the paper with discussions on how our technique differs from a previously-known strategy based on statistical inference.
\end{abstract}


\section{Introduction}
Anonymous (or `privacy-preserving') networks have been studied  in the computer science community, within the study of distributed and parallel computation,
since the late $1980$s \cite{YamashitaKameda_anonymous88}
and they have gained attention in the  control community in recent years 
due to the increasing importance of self-organized leaderless networks of sensors and of mobile robots (see the recent papers \cite{HendrickxEtAl-anonymous, VaragnoloEtAl_TACsubm} for a thorough summary of the literature
on anonymous networks from both communities). The peculiar aspect of anonymous networks, is that its agents are not able or not willing to provide an ID label uniquely
describing their identity. This might be due to technological limitations (when networks are self-configured, the agents do not run a time-consuming preliminary round to establish IDs, and it is impossible to ensure a priori that all devices of a given kind produced in the world have a unique identifier which they can use when self-configuring a network), or to privacy concerns (e.g., where agents represents smart-phones or computers, and the associated human being does not wish to reveal its participation within some network peer-to-peer activity).

Following the rich literature on consensus and on other distributed (i.e., leaderless, peer-to-peer) algorithms in sensor networks and in mobile robotic networks, we consider agents which are endowed with (limited) communication and computation capabilities, and have little or no knowledge of the network structure. In particular, the agents are able to send messages to their neighbors according to a given communication graph (which might depend on the inter-distances and possibly on obstacles, in the case of wireless communication).

In this setting, even simple tasks such as computing the number of nodes become very difficult. Negative results have been proved (see \cite{HendrickxEtAl-anonymous}), showing that in the anonymous setup it is impossible to have an algorithm able to compute the size of the network with probability one by using bounded communication, memory and computation-time resources. Then, the attention can be focused on finding algorithms able to compute the number of nodes with some good (albeit non-zero) probability, and/or to compute it using  resources growing with the size.

In \cite{VaragnoloEtAl_TACsubm,VaragnoloEtAl_CDC2010}, Varagnolo et al.\ propose a class of algorithms based on statistical inference and on well-known anonymous distributed algorithms for computing symmetric functions such as average and maximum (average-consensus and max-consensus). Their algorithms provide an estimate of the network size which is asymptotically unbiased, and whose variance can be made arbitrarily small at the price of increasing the communication, memory and computation resources.

In this paper, we propose a different technique to solve the problem of counting the number of nodes in the network. In our case, we use techniques from system identification: we ask the agents to run a linear dynamical system (suitably defined so that it can be constructed and run in the anonymous leaderless setup we are considering), and then to identify the order of the system, which will be a lower bound for the number of agents, and which will give the exact answer under some assumptions.

For numerical reasons that will be explained later,  the proposed algorithm is more effective when it is implemented using a finite field, instead of real numbers. Linear dynamical systems over finite fields have been studied for decades: most of the classical results have been stated for general fields in the book \cite{Kalman-book}, and a wide literature has appeared in the seventies and eighties related to the convolutional codes for error-correction in communications; in this paper we will use results from the recent paper \cite{sundaram2013structural} , to which we point the reader also for more references to classical literature.

The anonymous-network setup considered in this paper is described more in detail in the following subsection, while the proposed algorithm is introduced in Sect.~\ref{sec:main-idea-algo} and the version using finite fileds is discussed
in Sect.~\ref{sec:algo-galois}. In Sect.~\ref{sec:implement} we discuss some implementation details, and we provide and comment some simulation results.

\subsection{Problem setup}
We consider a network of $n$ agents. The communication constraint is described by a (directed) graph $\G$: agents are able to exchange messages with their neighbors in $\G$.
We assume that $\G$ is strongly connected, namely that for any pair of vertices $u,v$, there exists a directed path from $u$ to $v$. Moreover, we  assume that every vertex has a self-loop, which describes the fact that an agent always knows its own state.

Agents do not have any global knowledge about the network. However, we assume that
they know their neighborhood: they know the number of their in-neighbors, and moreover
they are able to assign labels to the edges in their neighborhood
(here we will assume they can assign labels to the incoming edges, but it is very simple to adapt
the algorithm to the case where they can assign label to their outgoing edges, which is the setup
considered in \cite{HendrickxEtAl-anonymous}).
This assumption requires
some minimal local coordination and identification of the neighbors, which is realistic in many
cases, but excludes some applications where there is a sharp privacy constraint and nodes
are required to process all incoming messages all outgoing messages irrespective of the sender (resp. receiver). This edge-labeling assumption is taken in \cite{HendrickxEtAl-anonymous}, and is 
required in part of the work in \cite{VaragnoloEtAl_CDC2010, VaragnoloEtAl_TACsubm},
since average-consensus algorithms require it (unless the graph is known to be balanced), while max-consensus
algorithms do not need it. An adaptation of our algorithm to the completely anonymous setup can be foreseen,
but its study is beyond the scope of this paper.

Our setup differs from \cite{HendrickxEtAl-anonymous}, since we allow
the nodes to have a memory which grows with the size of the network.
This apparently contradictive assumption (we cannot allocate a memory of size $n$ when we do not known $n$) 
in practice simply translates in the fact that memory limitations of the agents will be a key point to take into account
when evaluating the applicability of our algorithm to a given technology.

We assume that the graph is invariant during the duration of the algorithm, and that
there is some synchronization: although transmissions might not be simultaneous, 
each node knows the number of its neighbors, and is able to wait for all received messages
from iteration $k$ before performing iteration $k+1$. 
This excludes a gossip implementation.

\subsubsection*{Notation}
Throughout the paper, vectors will be column vectors and will be denoted with boldface lowercase fonts; for matrices we will use uppercase fonts. The symbol $\vt e_r$ denotes a column vector
with a $1$ in position $r$ and zero elsewhere, while $\ones$ denotes
a vector with all entries equal to~$1$.

\section{Counting agents as order-identification of a LTI system}
\label{sec:main-idea-algo}

\subsection{The main idea}

The main fact underlying this algorithm is that agents can run a discrete-time linear time-invariant system, in the same way used in the well-known linear consensus algorithm.
Each agent has a scalar state (denote by $x_u(k)$ the state of agent $u$ at time $k$);
agents send their own state to their neighbors, and then make a state-update which
is a linear combination of the received states:
\[ x_u(k+1) = \sum_{v \to u} A_{uv} x_v(k) \]
By stacking all $x_u(k)$'s in a vector $\vt x(k)$, this is a LTI system
$\vt x(k) = A \vt x(k)$,
where the state space has dimension $n$ (the number of agents),
and $A$ is a matrix consistent with the graph $\G$, i.e., $A_{uv}=0$
whenever $(v,u)$ is not an edge of $\G$.

The fact that each agent knows its own state, can be interpreted as saying
that agent $r$ can see the following output of the system:
\[  y^{(r)}(k) = \vt e_r^T \vt x(k) \,.\]
From the point of view of node $r$, it is then possible to construct
the Hankel matrix from the output data, as follows.
At time $k =2j-2$, using the outputs 
$y^{(r)}(0), \dots, y^{(r)}(2j-2)$,
node $r$ can construct the following $j \times j$ matrix:
\begin{equation} \label{eq:Hankel-definition}
 H^{(r)}_j =
\begin{bmatrix}
y^{(r)}(0)    & y^{(r)}(1) & \dots & y^{(r)}(j-1)\\
y^{(r)}(1)    & y^{(r)}(2) & \dots & y^{(r)}(j)    \\
\vdots &     &       & \vdots \\
y^{(r)}(j-1)& y^{(r)}(j) & \dots & y^{(r)}(2j-2)
\end{bmatrix} 
\end{equation}
The following proposition holds true; its simple proof is deferred to the appendix.
\begin{proposition} \label{prop:non-sing_leading_submatrices}
There exists a positive integer $\tilde n_r \le n$ such that
$H_j$ is invertible for all $j \le \tilde n_r$ and $H_j$ is singular
for all $j \ge \tilde n_r +1$.
$\hfill \square$ 
\end{proposition}
This means that node $r$ can compute $\tilde n_r$ by looking for $\tilde n_r+1$  the smallest
$j$ such that $H^{(r)}_j$ is singular.
Then, $\tilde n_r$ is always a lower bound for $n$.
The design of the matrix $A$ should be made in such a way to ensure that
$\tilde n_r = n$, as it is discussed in the next subsection.

\subsection{Decentralized design of the state-update matrix $A$}
The goal of the design is to choose a matrix $A$ such that
$H_j$ is invertible for all $j \le n$.
However, the design of $A$ should be done in a local way, without any global
knowledge of the graph, since nodes do not know the graph size and even the less they know
the overall structure of connections.
Agents are allowed to use only information about their neighbors.
Useful tools to achieve this goal come from the theory of structural controllability and observability of linear systems (see \cite{DionEtAl-survey-structural} for a survey of this research area), which suggests to use a simple
local random construction of the non-zero coefficients of $A$, achieving almost surely
$\tilde n_r = n $ for all $r$. The technique simply requires to sample
at random the coefficients of $A$, one in correspondence of each directed edge;
the random variables can be chosen from an arbitrary continuous distribution.

The following result holds.
\begin{theorem} \label{thm:structural}
Let $\G$ be strongly connected and with a self-loop on each vertex.
Let $m$ and $n$ be the number of its edges and nodes, respectively.
Let $A$ be a matrix consistent with $\G$ (i.e., $A_{uv} =0$ if $(v,u)$ is not an edge in $\G$), and let the $m$ elements of $A$ corresponding to the edges of $\G$ be free real-valued parameters
 $\lambda_1, \dots, \lambda_m$.
Let the initial state be $\vt x(0) = (\xi_1, \dots, \xi_n)^T$ for some free real-valued
parameters $\xi_1, \dots, \xi_n$.
Then,
$\rank H^{(r)}_n = n$ except for a set of parameters having zero Lebesgue measure in 
$\R^{m+n}$.~$\hfill \square$
\end{theorem}
The proof is analogous to the proofs of structural controllability and structural observability results (see \cite{DionEtAl-survey-structural} and the references therein). It is reported in the appendix.

\subsection{Proposed algorithm}

The considerations presented in the previous subsections lead us to proposing the following
algorithm.
\begin{algorithm}
\caption{Node-counting algorithm}\label{algo:node-counting}
\begin{algorithmic}[0]
\Statex  (Constructing A)
\State Each node $u$ sets a random label $A_{uv}$ on each incoming edge $(v,u)$
\Statex  (Initializing $\vt x(0)$)
\State Each node $u$ sets a random initial state $x_u(0)$
\Statex (Iterations $\vt x(k+1) = A \vt x(k) $  and Hankel matrix)
\For{$k=0,1, \dots$}
Each node $u$:
\State  sends its state $x_u(k)$ to its neighbors
\State  receives the state of its neighbors
\State  updates its own state:
		\[ x_u(k+1) = \sum_{v \to u} A_{uv} x_v(k) \]
\State  considers its state as an output $y^{(u)}(k) = x_u(k) $
	\If{$k$ is even and $\tilde n_r$ is undefined}
	 \State it constructs $H^{(u)}_{\frac{k}{2}+1}$ defined in (\ref{eq:Hankel-definition})
	 \If{$H^{(u)}_{\frac{k}{2}+1}$ is singular}
	 \State it returns $\tilde n_r = \frac{k}{2}$
	 \EndIf
	\EndIf 
\EndFor 
\end{algorithmic}
\end{algorithm}

Some remarks about this algorithm:
\begin{itemize}
\item The random variables can have any continuous distribution, and they 
aren't required to be independent, as far as the joint distribution remains continuous
\item 
There is no need to explicitly construct the Hankel matrix, which would require a memory space
of order $n^2$, see Sect.~\ref{sec:implement} for more detail on how to find the first singular Hankel matrix with a memory of size linear in $n$
\item 
The iterations for $k=0,1,\dots$ are not an infinite number of iterations leading asymptotically to the result: node $u$ can stop as soon as it finds a singular Hankel matrix, so that it stops after 
$\tilde n_r +1 \le n+1$ iterations. However, this detail has not been explicitly written in the algorithm, in order to leave the designer free to choose possible variations. For example, a node $u$, after having found a singular $H^{(r)}_{\frac{k}{2}+1}$, might stop testing the Hankel matrix singularity, but keep updating and broadcasting the state $x_u(k)$, in case some other node hasn't yet found a singular Hankel matrix. In this case, some other stopping criterion should be devised.
\end{itemize}

Theorem~\ref{thm:structural} guarantees that, with probability one, $\tilde n_r =n$ for all $n$,
i.e., all agents can correctly compute the size of the network.
Unfortunately, this theoretical result translates in a practically applicable algorithm only
for small size $n$, due to numerical issues. In fact,  when $n$ is larger than a few tens,
despite the theoretical guarantees from Proposition~\ref{prop:non-sing_leading_submatrices} and Theorem~\ref{thm:structural}, it often happens that numerically the Hankel matrix $H_j$ appears singular for $j$ much smaller than the correct value $n+1$.
This seems to be related to the fact that we are trying to identify a system of large size
from a scalar measurement, and therefore we need to let the system run for a long time,
so that the effect of the stability or instability of $A$ creates numerical problems.
In fact, if $A$ is unstable, then the entries of $\vt x(k)$ become very large, while if $A$ is stable they become very small. Even the choice of a marginally stable $A$ (such as the
stochastic matrices used in consensus algorithms, which ensure that all entries
of $\vt x(k)$ remain in the convex hull of the initial entries) does not help, since
the effect of the smaller stable eigenvalues of $A$ vanishes rapidly and is invisible at large $k$;
it is also intuitive that a system converging to consensus will produce for large k states $x_i(k)$
almost all equal, and with very little variation at next step $k+1$, thus making the Hankel matrix
having the last columns almost equal.
All such numerical issues do not appear in the case where Algorithm~\ref{algo:node-counting} is done performing all operations in some finite field, instead of in the field of real numbers, as it is discussed in the next section.

\section{Counting agents as order-identification of a LTI system over a finite field}
\label{sec:algo-galois}
The need to avoid the numerical problems arising when running the proposed algorithm
over the reals, suggests to choose some finite field instead, where there are no issues
about the effect of stability or instability of the system. An additional benefit is that the messages to be transmitted to neighbors can be exactly transmitted with a finite number of bits, and do not need approximations.

As a reminder, most results of linear algebra and of the theory of linear dynamical systems
are still true also over finite fields. However, special care should be taken for those
results whose proofs involve orthogonality or eigenvalues, which may fail, since the usual scalar product does not lead any more to a Hilbert space, and since it is no longer true that the
characteristic polynomial of a $n \times n$ matrix has $n$ solutions, which was true for complex numbers.
See \cite{sundaram2013structural} for all the results, in particular about controllability, which
are needed in our case.

Algorithm~\ref{algo:node-counting} is well-defined also in the case
where all variables and all parameters belong to a given finite field $\F_q$,
as far as all operations are done as defined in $\F_q$.
For the random parameters (entries of $A$ and of $\vt x(0)$),
they should be independent random variables, uniform on the field, since
under this assumption it is possible to guarantee that
$n$ is correctly computed with non-zero probability (Theorem~\ref{thm:structural-galois} below).

\begin{theorem} \label{thm:structural-galois}
Let $\G$ be strongly connected and with a self-loop on each vertex.
Let $m$ and $n$ be the number of its edges and nodes, respectively.
Let $A$ be a matrix consistent with $\G$ (i.e., $A_{uv} =0$ if $(v,u)$ is not an edge in $\G$), and let the $m$ elements of $A$ corresponding to the edges of $\G$ be free parameters
 $\lambda_1, \dots, \lambda_m$.
Let the initial state be $\vt x(0) = (\xi_1, \dots, \xi_n)^T$ for some free parameters $\xi_1, \dots, \xi_n$. Let $\lambda_1, \dots, \lambda_m, \xi_1, \dots, \xi_m$ be independent random
variables, uniformly distributed over the finite field $\F_q$
If $q \ge n^2$, then, for any given agent $r$
\[ \prob\Big(\rank H^{(r)}_n = n\Big) \ge 1 - \frac{n^2}{q} \,. \]
Moreover, defining $d = \frac{n^3+n}{2}$, if $q \ge d$, then
\[ \prob\Big(\rank H^{(r)}_n = n ~\forall r=1, \dots, n\Big) \ge 1 - \frac{d}{q} \,. \]
$\hfill \square$
\end{theorem}
The proof is inspired by the proofs of structural controllability over finite fields in \cite{sundaram2013structural}. It is reported in the appendix.

Notice that, differently from the real-valued case, it is not possible here to guarantee
a correct result with probability one. Moreover, the bound requires a field size which is very large for medium-size of $n$, e.g. when $n$ is of a few hundreds. However, simulation results show that the bound is very conservative, and it is possible to obtain probability of success larger than a half with field size around $2n$, as it is illustrated in Section~\ref{sec:simulations}.

\section{Algorithm implementation}
\label{sec:implement}

It is possible to test if an Hankel matrix $H_j$ constructed from
a string of values $y_0, y_1, \dots, y_{2j-2}$  is singular
without constructing and storing the whole matrix, and using only a memory linear
in $n$, instead of quadratic.  Algorithm~\ref{algo:Hankel-Lanczos}
is a slight modification of the Hankel-Lanczos factorization algorithm `AsymHankel' from \cite{Boley-Lee-Luk_LanczosAlgo},
where we have underlined the possibility  to re-use the computations
done for $H_j$ when testing $H_{j+1}$.
The algorithm works under the assumption that there exists a positive integer $\tilde n$
such that $H_{j}$ is invertible for all $j \le \tilde n$ and $H_{\tilde n} = 0$
(which is true for the Hankel matrices we are considering, see Prop.~\ref{prop:non-sing_leading_submatrices}). It finds such a $\tilde n$, because it constructs
$c_{1,1}, c_{2,2}, \dots, c_{\tilde n}, c_{\tilde n+1} = 0$ such that
$\det H_j = \prod_{i \le j} c_{i,i}$.
It can be successfully integrated with Algorithm~\ref{algo:node-counting},
since is processes the values $y_0, y_1, \dots$ as soon as they become available, and 
tests if $H_j$ is invertible as soon as $y_{2j-2}$ has been received.\\
\begin{algorithm}[!ht]
\caption{Finding the first singular Hankel matrix, 
based on (AsymHankel, \cite{Boley-Lee-Luk_LanczosAlgo})}
\label{algo:Hankel-Lanczos}
\begin{algorithmic}[0] 
\Statex Initialization
\State $c_{1,1}  = y_0$ 
\State $j = 1$
\Statex Iterations on $j$ until $c_{j,j} =0$ (i.e., until $\det H_j=0$) 
\While{$c_{j,j}\not=0$} 
\State (get $y_{2j-1}$)
\State $c_{2j,1} = y_{2j-1}$
	\If{$j>1$}
	\State $c_{2j-1,2} = c_{2j,1}-\gamma_{1,1} c_{2j-1,1}$
	\EndIf
	\If{$j>2$}
	\For{$k=2, \dots, j-1$}  
	\State $c_{2j-k, k+1} = $ 
			\State \;\;\; $c_{2j-k+1, k} 
				- \gamma_{k-1,k}c_{2j-k, k-1} 
				- \gamma_{k,k}c_{2j-k, k} $
	\EndFor 
	\EndIf
\State (get $y_{2j}$)
\State $c_{2j+1,1} = y_{2j}$
	\If{$j=1$}
	\State $\gamma_{1,1} = \frac{c_{2,1}}{c_{1,1}}$
	\Else  
	\State $\gamma_{j-1,j} =  \frac{c_{j,j}}{c_{j-1,j-1}}$
	\State $\gamma_{j,j} = \frac{c_{j+1,j}}{c_{j,j}} - \frac{c_{j,j-1}}{c_{j-1,j-1}} $
	\EndIf 
\State $c_{2j, 2} = c_{2j+1, 1} - \gamma_{1,1}c_{2j, 1}$ 
	\If{$j>1$}
	\For{$k=2, \dots, j$}  
	\State $c_{2j-k+1, k+1} = $ 
		\State $c_{2j-k+2, k} 
				-\gamma_{k-1,k}c_{2j-k+1, k-1} 
				- \gamma_{k,k}c_{2j-k+1, k} $
	\EndFor 
	\EndIf 
\State $j = j+1$
\EndWhile
\State \textbf{return} $\tilde n = j-1$ \Comment{$H_j$ is the first singular matrix}
\end{algorithmic}
\end{algorithm}

To prove that indeed $\det H_s = \prod_{i \le s} c_{i,i}$ for any $s \le \tilde n+1$, we can 
construct a lower-triangular matrix $L_s$ by letting
 $L_{ij} = c_{i,j}$ for all $1 \le j \le i \le s$.
A careful look at the algorithm (or a look at the algorithm AsymHankel in \cite{Boley-Lee-Luk_LanczosAlgo}, where the same computations are presented in matrix form, making this fact more visible) shows that, for all $j$, the $j$th column of $L_s$ is equal to the $j$th column of $H_s$
plus a linear combination of the columns 1st to $(j-1)$th of $H_s$. This implies that $L_s$ and 
$H_s$ have the same determinant.\\

Notice that the memory needed is linear in $n$, since at iteration
$j$ one needs all previously-computed $\gamma$'s,
and the $c$'s computed at iteration $j-1$ (while all other $c$'s can be discarded).
The total complexity is of order $n^2$, because  each iteration $j$ requires
a number of computations linear in $j$.

Finally notice that this algorithm can be applied also in the case of finite fields.

\section{Simulations and discussions}
\label{sec:simulations}
In this section we present some examples of results obtained with our algorithm. In all these examples, we use a graph  which is a circle where each node has a link towards the two nearest neighbors on its right and the two nearest nodes on its left, plus a self-loop. Other simulations with different graphs have lead us to similar remarks, and a deeper analysis of the effects of the graph topology (e.g., small-world, scale-free, Erdos-Renyi networks) on the performance is left for future work.

Moreover, for simplicity, we focus on finite fields whose size is a prime number $p$,
so that addition and multiplication are simply defined by operations modulo $p$.

Our simulations show that the bounds in Thm.~\ref{thm:structural-galois}
are very conservative. The theoretical result ensures that
a given node $r$ will correctly compute $\tilde n =n$ with probability at least $0.5$
if the size $q$ of the field is larger than $n^2$.
In our simulations, a success rate larger than $0.5$ is achieved when the field size is around $2n$.
Table~\ref{table:simul-galois-p2n} shows, for various values of $n$
and corresponding field size $p$ near to $2n$ (not exactly $2n$ since we require $p$ to be prime),
the success rate $S$, defined as the fraction of successful computations of $\tilde n_r = n$
from a given node $r$, in $1000$ realizations of the probabilistic algorithm.
Moreover, the success probability can be increased if the algorithm is modified as follows:
the algorithm is run twice (or, more in general, $M$ times) in parallel, with different realizations of $\vt x(0)$ and $A$, and then each node takes the maximum of the so-obtained bounds $\tilde n_r$.
The success rate of such a strategy is shown in Table~\ref{table:simul-galois-p2n}: $S_2$
is the fraction of pairs of realizations where at least one of the two gives $\tilde n_r = n$, over $1000$ pairs of realizations.
\begin{table}
\caption{For some values $n$ of the network size, and $p$ of the field size,
success rate $S$ = fraction of realizations where $\tilde n_1 = n$ and
$S_2$ = fraction of pairs of realizations where at least one of the two has $\tilde n_1 = n$.}
\label{table:simul-galois-p2n}
\begin{center}
\begin{tabular}{c|cccccccc|}
\hline
\!$n$\!      & 20  \!& 40  \!& 80  \!& 100 \!& 150 \!& 200 \!& 300 \!& 400\\
\hline
\!$p$\!      & 41  \!& 83  \!& 163 \!& 199 \!& 307 \!& 401 \!& 601 \!& 797\\
\hline
\!$S$\!      &0.582\!&0.587\!&0.618\!&0.608\!&0.649\!&0.611\!&0.580\!&0.593\\
\hline
\!\!$S_2$\!\!&0.834\!&0.859\!&0.848\!&0.837\!&0.856\!&0.844\!&0.834\!&0.841\\
\hline
\end{tabular}
\end{center}
\end{table}

\begin{figure}[th]
      \centering
      \includegraphics[width=.5\textwidth]{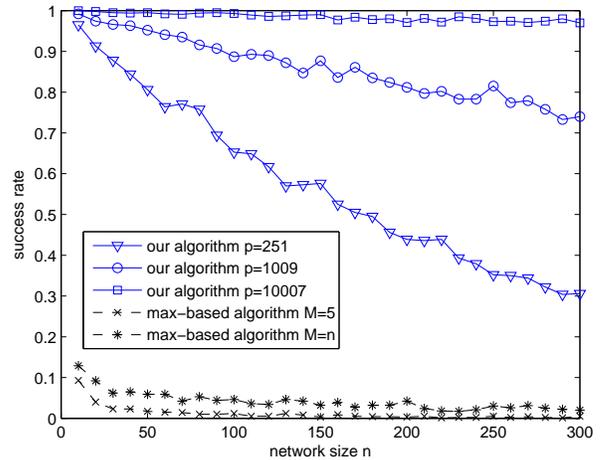}
      \caption{Success rate, i.e.\ fraction of realizations in which a given node $r$ has correctly computed the network size, in $1000$ realizations of each algorithm.}
      \label{fig:success-prob}
\end{figure}

\begin{figure}[th]
      \centering
      \includegraphics[width=.5\textwidth]{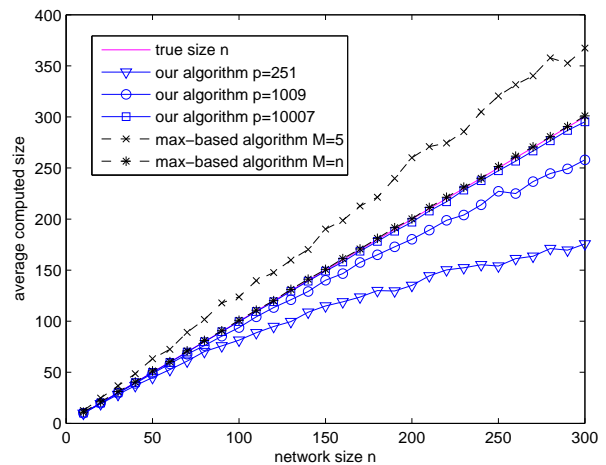}
      \caption{Average computed network size, over $1000$ realizations of each algorithm.}
      \label{fig:average-n}
\end{figure}

\begin{figure}[th]
      \centering
      \includegraphics[width=.5\textwidth]{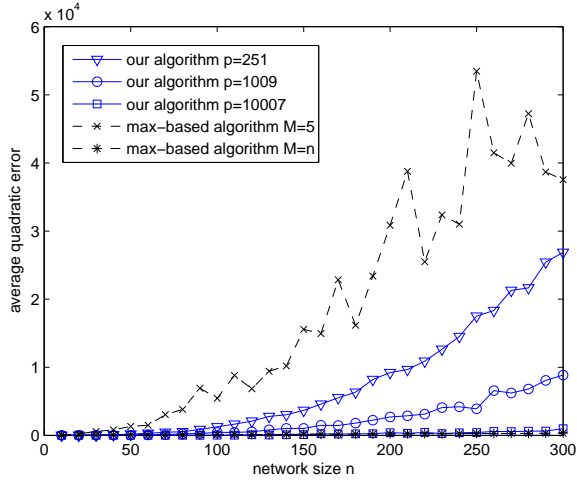}
      \caption{Average quadratic error, over $1000$ realizations of each algorithm.}
      \label{fig:quadratic-error}
\end{figure}

Another set of simulations showing that our algorithm performs well, and better than
predicted by the conservative bounds in Thm.~\ref{thm:structural-galois} is the following.
For a fixed field size $p$, we consider increasing values of $n$,
for each of which we run 1000 times our algorithm.
The solid lines with circles in Figures~\ref{fig:success-prob}, \ref{fig:average-n} and \ref{fig:quadratic-error}
depict the success rate (i.e., the fraction of realizations in which a given node $r$ has correctly computed $\tilde n_r =n$), the average computed value $\tilde n_r$, and then the average quadratic error (defined as the sum over all realizations of the square error $(\tilde n_r - n)^2$, divided by the number of realizations), respectively, 
for $p=251$, $p=1009$, and  $p = 10007$.
It is interesting to notice that even with $p=256 < n =300$ the success rate is  non-zero.


The same figures also show a comparison with one of the algorithms in \cite{VaragnoloEtAl_TACsubm}, which is based on statistical inference and on max-consensus (we have chosen this one, rather than the one based on average-consensus described in the same paper, because the authors prove that it has smaller variance). Such algorithm can be summarized as follows:
\begin{itemize}
\item Let $M$ be a positive integer, which indicates how many messages each node sends to its neighbors at each iteration.
\item Initialization:
each node $u$ extracts $M$ random variables $y_{1,u}, \dots, y_{M,u}$,
i.i.d. uniform in the interval $[0,1]$.
\item Max-consensus: a simple distributed anonymous algorithm allows
all agents to compute $f_1, \dots, f_M$,
where $f_i = \max_{u} y_{i,u}$.
\item  Given $f_1, \dots, f_M$, the maximum-likelihood estimate of $n$ is computed
as follows: 
\[ n_{\mathrm{ML}} = \round \left( \frac{M}{ -\sum_{i=1}^M \log f_i}  \right) \,, \]
where  $\round(\cdot)$ denotes rounding to the nearest integer.
\end{itemize}
The max-consensus algorithm consists in sending the current state to all neighbors,
and then updating the state to be the maximum of the received messages; this converges
to the exact maximum in $n$ iterations.
Hence, both our algorithm and the above one need a number of iterations linear in $n$,
but the latter does not provide the possibility to stop when $n$ iterations are done,
since $n$ is unknown, and thus it requires an upper bound of $n$ as a stopping criterion.
For memory requirements, our algorithm needs a memory linear in $n$ (roughly $6n$ numbers belonging to a field of size $q$), while the max-based algorithm requires a memory linear in $M$, where an increased $M$ results in a lower variance of the estimate.
For transmissions: at each iteration, each node sends only one value in our algorithm, and
$M$ values in the max-based one. We can also implement a version of our algorithm
where $M$ realizations are run in parallel, and hence $M$ messages are sent at each round (as it was done in Table~\ref{table:simul-galois-p2n} for $M=2$); this results in improved performance since the maximum of all computed sizes will give the correct size with higher probability, but it also incurs in increased memory, which is linear in the product $Mn$.
 
These considerations show that there is no natural obvious way to make a fair comparison between our algorithm and the max-based one, since depending whether memory or transmission complexity is more crucial, different choices can be made.
In Figures~\ref{fig:success-prob}, \ref{fig:average-n} and \ref{fig:quadratic-error}
we plot the success rate, the average computed network size, and the average quadratic error for our algorithm, and for the max-based algorithm with different values of $M$:  $M=5$ and $M=n$. The choice of $M=5$ has a transmission complexity comparable with our algorithm ($5$ times larger) and much smaller memory requirement; it gives very poor performance. The latter choice $M=n$ has a comparable memory requirement as our algorithm, and much higher transmission complexity ($n$ times larger); its performance is much worse in terms of probability of exact computation of $n$,  slightly better in terms of average $n$, and much better in terms of average quadratic error.

\appendix[\textbf{Proofs}]
Here we collect some proofs about properties (and structural properties) of the Hankel matrix.

First of all, notice that the following decomposition holds true:
\begin{equation} \label{eq:H=OC}
H^{(r)}_j = O^{(r)}_j C_j 
\end{equation}
where $O^{(r)}_j$ is the following $j \times n$ observability matrix
\[O^{(r)}_j = 
\begin{bmatrix}
\vt e_r^T \\
\hline \\
\vt e_r^T A\\
 \hline\\
\dots \\
\hline \\
\vt e_r^T A^{j-1}
\end{bmatrix}
\]
and $C_j$ is the following $n \times j$ matrix
\[ C_j =
\left[\begin{array}{c|c|c|c}
\vt x(0) & A \vt x(0) & \dots & A^{j-1} \vt x(0)
\end{array}\right]
\]
Notice that $C_j$ would be the controllability matrix if the system was $\vt x(k+1) = A \vt x(k) + B u(k)$, with $B = \vt x(0)$ and $u(k)$ a scalar input.

This decomposition means that $H^{(r)}_n$ is invertible if and only if the linear system $(A,\vt x(0), \vt e_r^T)$ is minimal (observable and controllable).
Moreover, this decomposition is useful in many proofs.

\textbf{Proof of  Prop.~\ref{prop:non-sing_leading_submatrices}}:
If $H_j$ is singular for some $j$, then at least one of the two matrices
$O^{(r)}_j $ and $C_j $ has rank smaller than $j$,
thus implying that at least one of $O^{(r)}_{j+1} $ and $C_{j+1} $ has rank smaller than $j+1$,
and hence also $H_{j+1}$ is singular.
The fact that $H_{n+1}$ is singular (and hence $\tilde n_r \le n$) follows from Cayley-Hamilton theorem, which ensures that $A^{n}$ is a linear combination of $I, A, A^2, \dots, A^{n-1}$,
so that the last line of $H_{n+1}$ is a linear combination of the previous ones.
$\hfill \blacksquare$\\

\textbf{Proof of  Theorems~\ref{thm:structural} and~\ref{thm:structural-galois}}:
This proof is inspired by the proofs of Thm.s 2 and 5 in \cite{sundaram2013structural}.
Notice that $\det H^{(r)}_n$ is a polynomial in the variables $\lambda_1, \dots, \lambda_m, \xi_1, \dots, \xi_n$.
The first part of the proof consists in proving that this is not the trivial all-zero polynomial.
To do so, we notice that $\det H^{(r)}_n = \det O^{(r)}_n \cdot \det C_n$,
where $\det O^{(r)}_n$ is a polynomial in the variables $\lambda_1, \dots, \lambda_m$,
which we will denote by $f_r(\lambda_1, \dots, \lambda_r)$,
and $\det C_n$ is a polynomial in the variables $\lambda_1, \dots, \lambda_m, \xi_1, \dots, \xi_n$,
which we will denote by $g(\lambda_1, \dots, \lambda_m, \xi_1, \dots, \xi_n)$.
We will now prove that neither of the two is the all-zero polynomial, by finding for each of the two (separately) a particular choice of the variables giving a non-zero evaluation.

Let's start with the determinant of $C_n$.
Choose $\xi_s = 1$ for one particular $s \in \{1, \dots, n\}$ and $\xi_h =0$ for all other $h$, so that $ C_n$ becomes equal to the controllability
matrix from node $s$, namely $C^{(s)}_n = \left[ \vt e_s | A \vt e_s | \dots | A^{n-1} \vt e_s \right]$.
Since $\G$ is strongly connected, it contains a directed spanning tree rooted at node $s$%
\footnote{A directed spanning tree rooted at $u$ is a subgraph of $\G$ containing exactly $n-1$ edges, 
and containing a directed path from $u$ to any other vertex.}.
Choose $\lambda_j$'s in the following way: $A_{ij}=1$ if $(j,i)$ is an edge of the spanning tree;
$A_{i,j}=0$ if $i \ne j$ and $(j,i)$ is not an edge of the spanning tree;
$A_{1,1}, \dots, A_{n,n}$ are non-zero distinct elements (this is possible, since by assumption the field has cardinality larger than $n$).
Consider a permutation reordering the vertices by non-decreasing depth (distance from vertex $s$) in the spanning tree. Notice that
$A = P^{-1} \tilde A P$, where $P$ is a permutation matrix corresponding to the above-mentioned re-ordering (so that, in particular,
$P \vt e_s = \vt e_1$ and $P^{-1} = P^T$), and $\tilde A$ is a lower-triangular matrix (since all its non-zero terms have $i \ge j$), with diagonal elements which are the same as the diagonal elements of $A$ (up to re-ordering) and hence are non-zero and all distinct.
Now
\begin{align*} 
C^{(s)}_n 
&= \left[  P^{-1} P \vt e_s |  P^{-1} \tilde A P \vt e_s | \dots | P^{-1} \tilde A^{n-1} P \vt e_s \right] \\
&= P^{-1} \left[   \vt e_1 | \tilde A  \vt e_1 | \dots | \tilde A^{n-1} \vt e_1 \right] 
\end{align*}
Since $\tilde A$ is a lower-triangular matrix  with distinct eigenvalues (the distinct diagonal elements), we can write $\tilde A = V^{-1} \Lambda V$, where $\Lambda = \diag (\tilde A) = \diag A$, and where $V$ is a lower-triangular matrix whose rows are left eigenvectors of $\tilde A$. Notice that we can choose  $V$ such that its first column is all-ones, so that $V \vt e_1 = \ones $, and hence
\begin{align*}
C^{(s)}_n 
&=
P^{-1} \left[  V^{-1} V \vt e_1 | V^{-1}\Lambda V  \vt e_1 | \dots | V^{-1}\Lambda^{n-1} V  \vt e_1 \right] \\
&=
P^{-1} V^{-1} M
\end{align*}
where $M$ is the Vandermonde matrix formed with the eigenvalues of $\tilde A$ (i.e., the diagonal elements of $A$):
\[ M =
[ \ones | \Lambda \ones | \dots | \Lambda^{n-1} \ones ]
=
\begin{bmatrix}
1 & A_{11} & A_{11}^2 & \dots & A_{11}^{n-1}\\
1 & A_{22} & A_{22}^2 & \dots & A_{22}^{n-1}\\
\vdots &   &          &       &   \vdots    \\
1 & A_{nn} & A_{nn}^2 & \dots & A_{nn}^{n-1}\\
\end{bmatrix}
\] 
Since the diagonal elements of $A$ are all distinct and non-zero, this Vandermonde matrix is invertible, and then also $C^{(s)}_n $ is invertible,
so that with the above-described choice of $\lambda_1, \dots, \lambda_m, x_1, \dots, x_n$, we have
$g(\lambda_1, \dots, \lambda_m, \xi_1, \dots, \xi_n) \ne 0$.
We have proved that $g$ is a non-zero polynomial. It is also useful to notice that its total degree is at most $\frac{n (n+1)}{2}$.
Indeed, each monomial in the calculation of the determinant of $C_n$ is the product of one entry per each column, 
and the entry from the $j$th column has degree at most $j-1$ in the variables $\lambda$'s and degree $1$ in the variables $\xi$'s,
so that the total degree is at most $0+1+ \dots + (n-1) = \frac{(n-1)n}{2}$ in the $\lambda$'s and $n$ in the $\xi$'s.\\

Now very similar considerations apply for the observability matrix $O^{(r)}_n$.
We construct the graph $\G^T$ defined by taking an edge $(i,j)$ if and only if $(j,i)$ is an edge of $\G$,
and we notice that also $\G^T$ is strongly connected and has self-loops at each vertex.
In this case, we consider a spanning tree in $\G^T$ rooted at node $r$, which corresponds to a set of paths
entering into node $r$ in $\G$. We make a similar construction as above, taking $A_{ij}$ to be 1 in correspondence with such paths,
with all-distinct diagonal elements, and zero elsewhere. Now, up to a permutation, $A$ is upper-triangular,
so that $A = P^{-1} \tilde A^T P$ with $\tilde A$ lower-triangular and $P$ a permutation matrix such that $P \vt e_r = \vt e_1$.
This also implies that $A^T = P^{-1} \tilde A P$ (since $P^T = P^{-1}$).
Now one can prove that  $\det (O^{(r)}_n)^T \ne 0$ exactly in the same way already used for $C^{(r)}_n$ above (replacing $A$ with $A^T$).
Hence, we have proved that  $f_r(\lambda_1, \dots, \lambda_m)$ is a non-zero polynomial. This polynomial has total degree at most $0+1+ \dots+ (n-1) = \frac{n(n-1)}{2}$.\\

Gathering the results on the two factors, we have finally proved that 
$\det H^{(r)}_n$ is a non-zero polynomial, with total degree at most $\frac{n(n-1)}{2} + \frac{n(n+1)}{2}  = n^2$.
This ends the first part of the proof.

Now, knowing that this isn't the trivial all-zero polynomial, we still need to prove
that the probability to find values of $\lambda_1, \dots, \lambda_m$ and $\xi_1, \dots, \xi_n$
that annihilate it is small enough.
In the case where the field is $\F = \R$, we can end the proof by noticing that the zeros of a polynomial
have zero Lesbesgue measure in  $\R^{m+n}$.
In the case where $\F = \F_q$ is a finite field of size $q \ge n^2$, the claim follows by the application of the following classic result.

\begin{lemma}[Schwartz-Zippel] Given a non-zero polynomial $p(z_1, \dots, z_{\ell}) \in \F_q[z_1, \dots, z_{\ell}]$, if its total degree is at most $d$, the field size is $q \ge d$ and $z_1, \dots, z_n$ are iid unif in $\F_q$ then
\[ \prob\big(p(z_1, \dots, z_{\ell}) = 0\big) \le \frac{d}{q} \,.\]
\end{lemma}

This ends the proof of the results concerning the probability that a given node $r$
correctly computes $\tilde n_r = n$.
When considering the probability that all agents simultaneously get the correct $n$,
just notice that the event `$\det H^{(r)}_n \ne 0$ for all $r = 1, \dots, n$' is equivalent to the following:
\[ \det O^{(1)}_n \cdot \det O^{(2)}_n \cdot \dots \det O^{(n)}_n \cdot\det C_n \ne 0 \,,\]
since a product is non-zero if and only if all factors are non-zero,
and since the factorization (\ref{eq:H=OC}) is true for all $r$ with the same
 $C_n$.

Then, the proof is the same as above, and the total degree of the polynomial 
$\det O^{(1)}_n \cdot \det O^{(2)}_n \cdot \dots \det O^{(n)}_n \cdot\det C_n $ is
at most $d = \frac{n(n-1)}{2} n + \frac{n(n+1)}{2}  = \frac{n^3+n}{2}$.

$\hfill \blacksquare$\\


\bibliographystyle{plain}
\bibliography{Garin-Yuan-bibtex}

\end{document}